\documentclass{article}

\usepackage{amssymb}
\usepackage{amssymb,amsmath,latexsym}
\begin{document}
\title{\bf \small{A new study on the absolute matrix summability of non-decreasing sequences}}
\author{\small \c{S}ebnem YILDIZ\\
\small Department of Mathematics,
Ahi Evran University, K{\i}r\c{s}ehir, Turkey\\
\small  e-mail: sebnemyildiz@ahievran.edu.tr; sebnem.yildiz82@gmail.com
\date{}} \maketitle
\date{}
\noindent {\bf Abstract.} Recently, in \cite{Bor4}, Bor proved a main theorem dealing with $|\bar{N}, p_{n}|_{k}$ summability of non-decreasing sequences. In the present paper, we have generalized that theorem for  $|A, p_{n}|_{k}$ summability method by using almost increasing sequences and taking normal matrices in place of weighted mean matrices.
\section {Introduction}
Let $\sum a_{n}$ be a given infinite series with partial sums $(s_{n})$. We denote by $u_{n}^{\alpha}$ the nth Ces\`aro mean of order $\alpha$, with $\alpha>-1$, of the sequence $(s_{n})$, that is (see \cite{Ce})
\begin{align}
u_{n}^{\alpha}=\frac{1}{A_{n}^{\alpha}}\sum_{v=0}^{n}A_{n-v}^{\alpha-1}s_{v}~~~~
\end{align}
\noindent where
\begin{align} A_{n}^{\alpha}=\frac{(\alpha+1)(\alpha+2)....(\alpha+n) }{n!}=
O(n^\alpha), \quad {A_{-n}^\alpha=0} \quad {\texttt{for}} \quad  n>0.\
\end{align}
A series $\sum{a_n}$ is said to be summable
$\mid{C},\alpha\mid_k$, $k\geq1$, if (see \cite{Fl})
\begin{align}
\sum_{n=1}^{\infty} n^{k-1}\mid
u_n^{\alpha}-u_{n-1}^{\alpha}\mid^{k}<\infty.
\end{align}
If we set $\alpha$=1, then we have ${\mid{C,1}\mid}_k$ summability. Let $(p_{n})$ be a sequence of positive number such that
\begin{align}
P_{n}=\sum_{v=0}^{\infty}p_{v}\rightarrow \infty \quad {as}\quad
{n}\rightarrow \infty ,\quad (P_{-i}=p_{-i}=0,~~i\geq1).
\end{align}
\small {\it{ 2010 AMS Subject Classification:} 26D15, 40D15, 40F05, 40G99, 42A24, 46A45}.\\
\small{\it{ Keywords:} \small Riesz mean, absolute matrix summability, summability factors, infinite series,
H\"older inequality, Minkowski inequality}.\\
The sequence-to-sequence transformation\\
\begin{align}
w_{n}=\frac{1}{P_{n}} \sum_{v=0}^{n}p_{v}s_{v}
\end{align}
defines the sequence $(w_{n})$ of the Riesz mean or simply the $\left(\bar{N},p_{n}\right)$
mean of the sequence $(s_{n})$, generated by the sequence of coefficients $(p_{n})$ (see \cite{Ha}).
The series $\sum{a_n}$ is said to be summable
$|\bar{N},p_{n}|_{k}$, $k\geq1$, if (see \cite{Bor1})
\begin{align}
\sum_{n=1}^{\infty}\left(\frac{P_{n}}{p_{n}}\right)^{k-1}\left|w_{n}-w_{n-1}\right|^{k} < \infty.
\end{align}
In the special case when $p_{n}=1$ for all values of $n$ (respect. $k=1$), then $|\bar{N},p_{n}|_{k}$ summability is the same as $|C,1|_{k}$ (respect. $|\bar{N},p_{n}|$) summability.\\
Let $A=(a_{nv})$ be a normal matrix,
i.e., a lower triangular matrix of nonzero diagonal entries. Then $A$ defines the sequence-to-sequence transformation, mapping the
sequence $s=(s_{n})$ to $As=\left(A_{n}(s)\right)$, where
\begin{align}\label{eq:6}
A_{n}(s)=\sum_{v=0}^{n}a_{nv}s_{v}, \quad n=0,1,...
\end{align}
The series $\sum a_{n}$ is said to be summable  $\left|A,p_{n}\right|_{k}$, $k\geq 1$,
if (see \cite{WT})
\begin{align}\label{eq:8}
\sum_{n=1}^{\infty}\left(\frac{P_{n}}{p_{n}}\right)^{k-1} \left|\bar{\Delta}A_{n}(s)\right|^{k}<\infty.
\end{align}
where
\begin{align}\label{eq:9}
\bar{\Delta}A_{n}(s)=A_{n}(s)-A_{n-1}(s).
\end{align}
 \noindent  Note that in the special case if we take $ p_{n}=1$ for all $n$, $\left|A,p_{n}\right|_{k}$ summability is the same as  $\left| A\right|_{k}$ summability (see \cite{NT}). Also,  if we take $a_{nv}=\frac{p_{v}}{P_{n}}$, then $\left| A,p_{n}\right|_{k}$ summability reduces to  $\left|\bar{N},p_n\right| _k$ summability. Furthermore, if we take $a_{nv}=\frac{p_{v}}{P_{n}}$ and $p_{n}=1$ for all values of $n$, then $\left|A,p_{n}\right|_{k}$ summability is the same as $\left| C,1\right|_k$ summability.
\section {The Known Results} 

A positive sequence $(b_{n})$ is said to be almost increasing if there exists a positive increasing sequence $(z_{n})$ and two positive constants 
$A$ and $B$ such that $Az_{n}\leq b_{n}\leq Bz_{n}$ (see \cite{bari}). It is known that every increasing sequences is an almost increasing sequence but the converse need not be true. 
Quite recently, Bor has  proved the following theorems concerning on summability factors of the absolute weighted mean.

{\bf Theorem 2.1} \cite{Bor2} Let $(X_{n})$ be a positive non-decreasing sequence and suppose that there exists sequences $(\beta_{n})$ and $(\lambda_{n})$  such that
\begin{align}
&|\Delta\lambda_{n}|\leq\beta_{n},\\
&\beta_{n}\rightarrow 0 \quad \text{as} \quad n\rightarrow \infty\\
&\sum_{n=1}^{\infty}n|\Delta\beta_{n}|X_{n}=O(1),\\
&|\lambda_{n}|X_{n}=O(1).
\end{align}
If
\begin{align}
\sum_{n=1}^{m}\frac{|s_{n}|^{k}}{n}=O(X_{m}) \quad \text{as} \quad m\rightarrow \infty,
\end{align}
and $(p_{n})$ is a sequence that
\begin{align}
&P_{n}=O(np_{n}),\\
&P_{n}\Delta p_{n}=O(p_{n}p_{n+1}),
\end{align}
then the series $\sum a_{n}\frac{P_{n}\lambda_{n}}{np_{n}}$ is summable $|\bar{N},p_{n}|_{k}$, $k\geq1$.\\
Later on, Bor  has recently proved the following theorem using under weaker conditions.\\
{\bf Theorem 2.2} \cite{Bor4} Let $(X_{n})$ be a positive non-decreasing sequence. If the sequences $(X_{n})$ , $(\beta_{n})$, $(\lambda_{n})$, and $(p_{n})$
satisfy the conditions (10)-(13), (15)-(16), and 
\begin{align}
\sum_{n=1}^{m}\frac{|s_{n}|^{k}}{nX_{n}^{k-1}}=O(X_{m}) \quad \text{as} \quad m\rightarrow \infty,
\end{align}
then the series $\sum a_{n}\frac{P_{n}\lambda_{n}}{np_{n}}$ is summable $|\bar{N},p_{n}|_{k}$, $k\geq1$.
 \section {The Main Results}

  The aim of this paper is to generalize Theorem 2.2 for $|A,p_{n}|_{k}$ summability factors using almost increasing sequences in place of positive non-decreasing sequence. So, we have generalized Theorem 2.2 under weaker hypothesis by using normal matrices.
 \\Given a normal matrix $A=(a_{nv})$, we associate two lower
 semimatrices $\bar{A}=(\bar{a}_{nv})$ and $\hat{A}=(\hat{a}_{nv})$
 as follows:
\begin{align}
\bar{a}_{nv}=\sum_{i=v}^{n}a_{ni},\quad n,v=0,1,...\quad \bar{\Delta}a_{nv}=a_{nv}-a_{n-1,v}, \quad a_{-1,0}=0
\end{align}
and
\begin{align}
\hat{a}_{00}=\bar{a}_{00}=a_{00},\quad
\hat{a}_{nv}=\bar{\Delta}\bar{a}_{nv},\quad n=1,2,...
\end{align}
 It may be noted that $\bar{A}$ and $\hat{A}$ are the well-known
 matrices of series-to-sequence and series-to-series
 transformations, respectively. Then, we have
 \begin{align}\label{eq:13}
 A_{n}(s)&=\sum_{v=0}^{n}a_{nv}s_{v}=
 \sum_{v=0}^{n}\bar{a}_{nv}a_{v}
 \end{align} 
 and
 \begin{align}
 \bar{\Delta}A_{n}(s)&=\sum_{v=0}^{n}\hat{a}_{nv}a_{v}.
 \end{align}
 With this notation we have the following theorem.\\
{\bf Theorem 3.1} Let  $A=(a_{nv})$ be a positive normal matrix such that
 	\begin{align}\label{eq:15}
 	\overline{a}_{n0}&=1,\     n=0,1,...,\\
 	a_{n-1,v}&\geq a_{nv},\ \textnormal{for}~~   n\geq v+1,\\
 	a_{nn}&=O(\frac{p_{n}}{P_{n}}),\\
 	\sum_{v=1}^{n-1}a_{vv}\hat{a}_{n,v+1}&=O(a_{nn})
 	\end{align}
and let $(X_{n})$ be an almost increasing sequence. If the sequences $(X_{n})$, $(\beta_{n})$, $(\lambda_{n})$, and $(p_{n})$
satisfy the conditions of Theorem 2.2, then the series $\sum a_{n}\frac{P_{n}\lambda_{n}}{np_{n}}$  is summable $\left|A,p_{n}\right|_{k}$, $k\geq 1$.\\
We need the following lemmas for the proof of Theorem 3.1.\\
{\bf Lemma 3.1} \cite{KN} Under the conditions on  $(X_{n})$, $(\beta_{n})$, and $(\lambda_{n})$ as expressed in the statement of Theorem 2.2, we have the following:
\begin{align}
nX_{n}\beta_{n}& =O(1),\\
\sum_{n=1}^{\infty}\beta_{n} X_n& <\infty.
\end{align}\\
{\bf Lemma 3.2} \cite{KN2} If the conditions (15) and (16) of Theorem 2.1 are satisfied, then $\Delta\left(\frac{P_{n}}{np_{n}} \right) =O\left( \frac{1}{n}\right).$ \\
{\bf Remark} Under the conditions on the sequence $(\lambda_{n})$ of Theorem 2.1, we have that $(\lambda_{n})$ is bounded and $\Delta\lambda_{n}=O(1/n)$ (see \cite{Bor2}).
\section{Proof of Theorem 3.1}
Let $(V_{n})$ denotes the A-transform of the series $\sum a_{n}\frac{P_{n}\lambda_{n}}{np_{n}}$. Then, by the definition, we have that
\begin{align*}
	\bar{\Delta}V_{n}& = \sum_{v=1}^{n}\hat{a}_{nv}a_{v}\frac{P_{v}\lambda_{v}}{vp_{v}}.
\end{align*}
Applying Abel's transformation to this sum, we have that
\begin{align*}
	\bar{\Delta}V_{n}& =\sum_{v=1}^{n-1}\Delta_{v}\left( \frac{\hat{a}_{nv}P_{v}\lambda_{v}}{vp_{v}}\right)\sum_{r=1}^{v}a_{r}+\frac{\hat{a}_{nn}P_{n}\lambda_{n}}{np_{n}}\sum_{r=1}^{n}a_{r}\\
	\bar{\Delta}V_{n}& =\sum_{v=1}^{n-1}\Delta_{v}\left( \frac{\hat{a}_{nv}P_{v}\lambda_{v}}{vp_{v}}\right)s_{v}+\frac{\hat{a}_{nn}P_{n}\lambda_{n}}{np_{n}}s_{n},	
\end{align*}
by the formula for the difference of products of sequences (see \cite{Ha}) we have
\begin{align*}
\bar{\Delta}V_{n}&=\frac{a_{nn}P_{n}\lambda_{n}}{np_{n}}s_{n}+\sum_{v=1}^{n-1}\frac{P_{v}\lambda_{v}}{vp_{v}}\Delta_{v}(\hat{a}_{nv})s_{v}+\sum_{v=1}^{n-1}\hat{a}_{n,v+1}\lambda_{v}\Delta\left( \frac{P_{v}}{vp_{v}}\right)s_{v}+\sum_{v=1}^{n-1}\hat{a}_{n,v+1}\frac{P_{v+1}}{(v+1)p_{v+1}}\Delta\lambda_{v}s_{v}\\
\bar{\Delta}V_{n}&=V_{n,1}+V_{n,2}+V_{n,3}+V_{n,4}.
\end{align*}
To complete the proof of Theorem 3.1, by Minkowski's inequality, it is sufficient to show that
\begin{align}
	\sum_{n=1}^{\infty}\left(\frac{P_{n}}{p_{n}}\right)^{k-1}\mid V_{n,r}\mid^{k}< \infty, \quad \textnormal{for} \quad{r=1,2,3,4.}
\end{align}
Firstly, by applying Abel's transformation and in view of the hypotheses of Lemma 3.1, Lemma 3.2, and Theorem 3.1 we complete the proof of Theorem 3.1 .\\
\section{Conclusions}
1. If we take $(X_{n})$ as a positive non-decreasing sequence and  $a_{nv}=\frac{p_{v}}{P_{n}}$ in Theorem 3.1, then we obtain Theorem 2.2 and if we put $k=1$ in Theorem 2.2, we have a known result of Mishra and Srivastava dealing with  $\left| \bar{N},p_{n}\right| $ summability factors of infinite series (see \cite{KN2}).\\
2. If we take $p_{n}=1$ for all values of $n$ in Theorem 3.1, then we get a new result dealing with the $|A|_{k}$ summability method.\\
3. If we take $a_{nv}=\frac{p_{v}}{P_{n}}$ and $p_n=1$ for all values of $n$ in Theorem 3.1, then we obtain  a known  result of Mishra and Srivastava  concerning the $\mid C,1\mid_k$ summability factors of infinite series (see \cite{KN1}).\\

\end{document}